\documentclass[submission]{dmtcs}
\usepackage[cp1251, cp866]{inputenc}
\usepackage[english]{babel}
\usepackage{amsmath, amssymb}
\usepackage{graphicx}
\usepackage{varioref}
 \newcommand{\pf}{\noindent {\it Proof:}  }

 \newtheorem{tm}{Theorem}

 \newtheorem{lemma}[tm]{Lemma}
 \newtheorem{cor}[tm]{Corollary}
 \newtheorem{prop}[tm]{Proposition}

 \renewcommand{\a}{\hspace{0.6cm}}
\renewcommand{\b}{\hspace{1.15cm}}
\renewcommand{\c}{\hspace{1.75cm}}
\renewcommand{\d}{\hspace{2.3cm}}

 \author{Olga Glebova\addressmark{1}{\ },
 Yury Metelsky \addressmark{1}{\ }
 \and Pavel Skums\addressmark{1}{\ }}
\title[Krausz dimension and its generalizations in special graph classes]{Krausz dimension and its generalizations in special graph classes}
\address{
\addressmark{1} Department of Mechanics and Mathematics,
 Belarus State University, 4 Nezavisimosti av.,
 220030 Minsk, Republic of Belarus.\\
 E-MAIL: glebovaov@gmail.com
}
 \revision{1} \received{???} \revised{???} \accepted{???}
\begin{document}
\maketitle
\begin{abstract}
A {\it krausz $(k,m)$-partition} of a graph $G$ is the partition
of $G$ into cliques, such that any vertex belongs to at most $k$
cliques and any two cliques have at most $m$ vertices in common.
The {\it $m$-krausz} dimension $kdim_m(G)$ of the graph $G$ is the
minimum number $k$ such that $G$ has a krausz $(k,m)$-partition.
1-krausz dimension is known and studied krausz dimension of graph
$kdim(G)$.

In this paper we prove, that the problem $"kdim(G)\leq 3"$ is
polynomially solvable for chordal graphs, thus partially solving
the problem of P. Hlineny and J. Kratochvil. We show, that the
problem of finding $m$-krausz dimension is NP-hard for every
$m\geq 1$, even if restricted to (1,2)-colorable graphs, but the
problem $"kdim_m(G)\leq k"$ is polynomially solvable for
$(\infty,1)$-polar graphs for every fixed $k,m\geq 1$.
\end{abstract}

\section{Introduction}
\label{sec:in}

In this paper we consider finite undirected graphs without loops
and multiple edges. The vertex and the edge sets of a graph
(hypergraph) $G$ are denoted by $V(G)$ and $E(G)$ respectively.
$N(v)=N_G(v)$ is the neighborhood of a vertex $v$ in $G$ and
$deg(v)$ is the degree of $v$. Let $G(X)$ denote the subgraph of
$G$ induced by a set $X\subseteq V(G)$ and $ecc_G(v)$ is the
eccentricity of a vertex $v\in V(G)$.

A {\it krausz partition} of a graph $G$ is the partition of $G$
into cliques (called {\it clusters} of the partition), such that
every edge of $G$ belongs to exactly one cluster. If every vertex
of $G$ belongs to at most $k$ clusters then the partition is
called {\it krausz $k$-partition}. The {\it krausz dimension}
$kdim(G)$ of the graph $G$ is a minimal $k$ such that $G$ has
krausz $k$-partition.

Krausz $k$-partitions are closely connected with the
representation of a graph as an intersection graph of a
hypergraph. The {\it intersection graph} $L(H)$ of a hypergraph
$H=(V(H),E(H))$ is defined as follows:

\begin{itemize}
\item[1)] the vertices of $L(H)$ are in a bijective correspondence
with the edges of $H;$

\item[2)] two vertices are adjacent in $L(H)$ if and only if the
corresponding edges have a nonempty intersection.
\end{itemize}

Hypergraph $H$ is called {\it linear}, if any two of its edges
have at most one common vertex; $k$-{\it uniform}, if every edge
contains $k$ vertices.

The {\it multiplicity} of the pair of vertices $u$, $v$ of the
hypergraph $H$ is the number $m(u,v)=|\{\mathcal{E}\in E(H) :
u,v\in \mathcal{E}\}|$; the {\it multiplicity} $m(H)$ of the
hypergraph $H$ is the maximum multiplicity of the pairs of its
vertices. So, linear hypergraphs are the hypergraphs with the
multiplicity 1.

Denote by $H^{*}$ the dual hypergraph of $H$ and by $H_{[2]}$ the
2-section of $H$ (i.e. the simple graph obtained from $H$ by
transformation each edge into a clique). It follows immediately
from the definition that

\begin{equation}\label{BergeThm}
    L(H)=(H^*)_{[2]}
\end{equation}

\noindent (first this relation was implicitly formulated by C.
Berge in \cite{B89}). This relation implies that a graph $G$ has
krausz $k$-partition if and only if it is intersection graph of
linear $k$-uniform hypergraph.

A graph is called $(p,q)$-{\it colorable} \cite{Bran96}, if its
vertex set could be partitioned into $p$ cliques and $q$ stable
sets. In this terms $(1,1)$-colorable graphs are well-known split
graphs.

Another generalization of split graphs are {\it polar graphs} (see
\cite{DEHS08},\cite{TC78}). A graph $G$ is called {\it polar} if
there exists a partition of its vertex set
\begin{equation}\label{DefinitionPolar}
    V(G)=A\cup B,\ A\cap B =\emptyset
\end{equation}
({\it bipartition} $(A,B)$) such that all connected components of
the graphs $\overline{G}(A)$ and $G(B)$ are complete graphs. If,
in addition, $\alpha$ and $\beta$ are fixed integers, and the
orders of connected components of the graphs $\overline{G}(A)$ and
$G(B)$ are at most $\alpha$ and $\beta$ respectively, then the
polar graph $G$ with bipartition (\ref{DefinitionPolar}) is called
$(\alpha,\beta)$-{\it polar}. Given a polar graph $G$ with
bipartition (\ref{DefinitionPolar}), if the order of connected
components of the graph $\overline{G}(A)$ (the graph $G(B)$) is
not restricted above, then the parameter $\alpha$ (respectively
$\beta$) is supposed to be equal $\infty$. Thus an arbitrary polar
graph is $(\infty,\infty)$-polar, and a split graph is
$(1,1)$-polar.

Denote by $KDIM(k)$ the problem of determining whether
$kdim(G)\leq k$ and by $KDIM$ the problem of finding the krausz
dimension.

The class of line graphs (intersection graphs of linear 2-uniform
hypergraphs, i.e graphs with krausz dimension at most 2) have been
studied for a long time. It is characterized by a finite list of
forbidden induced subgraphs \cite{Bein}, efficient algorithms for
recognizing it (i.e. solving the problem $KDIM(2)$) and
constructing the corresponding krausz 2-partition are also known
(see for example \cite{DS95}, \cite{L73}, \cite{NN90},
\cite{R73}).

 The situation changes radically if one takes $k=3$ instead
 of $k=2$ : the problem $KDIM(k)$
is NP-complete for every fixed $k\geq 3$ \cite{HK97}. The case
$k=3$ was studied in the different papers (see
\cite{JKL97},\cite{MT97},\cite{NRSS80},\cite{NRSS82},\cite{SST09}),
and several graph classes, where the problem $KDIM(3)$ is
polynomially solvable or remains NP-complete, were found.

In \cite{HK97} P.~Hlineny and J.~Kratochvil studied the
computational complexity of the krausz dimension in detail.
Besides another results, the following results were obtained in
their paper:

\begin{itemize}

\item[1)] The problem $KDIM$ is polynomially solvable for graphs
with bounded treewidth. In particular, it is polynomially solvable
for chordal graphs with bounded clique size.

\item[2)] For the whole class of chordal graphs the problem
$KDIM(k)$ is NP-complete for every $k\geq 6$.

\end{itemize}

 So, the problem of deciding the complexity of $KDIM(k)$
restricted to chordal graphs for $k=3,4,5$ was posed by P.~Hlineny
and J.~Kratochvil. As a partial answer to it, in the Section
\ref{sec:chorddim3} we prove that the problem $KDIM(3)$ is
polynomially solvable in the class of chordal graphs.

In the Section \ref{mkrausz} we consider the natural
generalization of the krausz dimension. The {\it krausz
$(k,m)$-partition} of a graph $G$ is the partition of $G$ into
cliques (called {\it clusters} of the partition), such that any
vertex belongs to at most $k$ clusters of the partition, and any
two clusters have at most $m$ vertices in common. As above, the
relation (\ref{BergeThm}) implies the following statement:
\begin{prop}\label{BergeThm+} A graph $G$ has krausz $(k,m)$-partition
if and only if it is the intersection graph of a $k$-uniform
hypergraph with the multiplicity at most $m$.
\end{prop}
The {\it $m$-krausz} dimension $kdim_m(G)$ of the graph $G$ is the
minimum $k$ such that $G$ has a krausz $(k,m)$-partition. The
krausz dimension in this terms is the 1-krausz dimension.

Denote by $KDIM_m$ the problem of determining the $m$-krausz
dimension of graph, by $KDIM_m(k)$ the problem of determining
whether $kdim_m(G)\leq k$ and by $L_k^m$ the class of graphs with
a krausz $(k,m)$-partition. It was proved in \cite{LevTysh93} that
the class $L_3^m$ could not be characterized by a finite set of
forbidden induced subgraphs for every $m\geq 2$, but the
complexity of the problem $KDIM_m$ for an arbitrary $m$ was not
established yet. We prove that the problem $KDIM_m$ is NP-hard for
every $m\geq 1$, even if restricted to the class of
$(1,2)$-colorable graphs.

The class $L_k^m$ is hereditary (i.e. closed with respect to
deleting the vertices) and therefore can be characterized in terms
of forbidden induced subgraphs. We prove that for every fixed
integers $m,k$ such finite characterization of the class exists
when restricted to $(\infty,1)$-polar graphs. In particular, it
follows that the problem $KDIM_m(k)$ is polynomially solvable for
$(\infty,1)$-polar graphs for every fixed $m$ and $k$. In
particular, it generalizes the result of \cite{HK97} and
\cite{M97}, that for every fixed $k$ the problem $KDIM(k)$ is
polynomially solvable for split graphs.

\section{Krausz 3-partitions of chordal graphs}
\label{sec:chorddim3}

Let $F$ be a family of cliques of graph $G.$ The cliques from $F$
are called {\it clusters of $F.$} Denote by $l_F(v)$ the number of
clusters from $F$ covering the vertex $v.$

A maximal clique with at least $k^2-k+2$ vertices is called a {\it
$k$-large} clique. For such cliques the following statement holds:

\begin{lemma}\label{large}\cite{HK97,JKL97,NRSS82} Any $k$-large
clique of a graph $G$ belongs to every krausz $k$-partition of
$G.$
\end{lemma}

Further in this section 3-large clique will be called simply {\it
large clique}.

Let $G$ be a graph with $kdim(G)\leq 3$ and $Q$ be some its krausz
3-partition. Any subset $F \subseteq Q$ is called {\it a fragment
of the krausz 3-partition} $Q$ (or simply {\it a
 fragment}).

Let $F$ be some fragment of krausz 3-partition $Q$ and $H$ be the
subgraph of $G$ obtained by deleting edges covered by $F$ ($F$
could be empty). Fix some vertex $a\in V(H)$ and positive integer
$k$. Denote by $B_k[a]$ the $k$th neighborhood of  $a$ in $H$,
i.e. the set of vertices at distance at most $k$ from $a$. A
family of cliques $F_k(a)$ in $H$ is called {\it $(a,k)$-local
fragment} (or simply a {\it local fragment}), if

\begin{itemize}
\item[(1)] every edge with at least one end in $B_k[a]$ is covered
by some cluster of $F_k(a)$;

\item[(2)] every vertex $v\in B_k[a]$ belongs to at most
$3-l_F(v)$ clusters of $F_k(a)$.

\item[(3)] every two clusters of $F_k(a)$ have at most one common
vertex.
\end{itemize}

A clique $C$ is called {\it special}, if $C$ is a cluster of every
$(a,k)$-local fragment for some $a$ and $k$. In particular, by
Lemma \ref{large} large cliques are special.

The following statements are evident.

\begin{lemma}\label{technical}

\begin{itemize}
\item[1)] If $deg(v)\geq 19$ for some vertex $v\in V(G)$, then $v$
is contained in some large clique.

\item[2)] If $l_F(v) = 2$, then $C = N_{H}(v)\cup\{v\}$ is a
special clique.

\item[3)] If $v\in B_k[a]$ is adjacent to at least $4-l_F(v)$
vertices of the cluster $C$ of some local fragment $F_k(a)$, then
$v\in C$.

\item[4)] for every $a\in V(H)$ and every $k$ there exists at
least one $(a,k)$-local fragment;

\item[5)] If the clique $C$ is special, then $F\cup \{C\}$ is a
fragment.

\end{itemize}

\end{lemma}

\pf Let's illustrate, for example, 3) and 5). If $v\in B_k[a]$ is
adjacent to vertices $v_1,...,v_{4-l_F(v)}\in C\in F_k(a)$, but
$v\not\in C$, then the edges $vv_1,...vv_{4-l_F(v)}$ should be
covered by different clusters of $F_k(a)$. It contradicts (2).

The family of cliques $X = \{C\in Q\setminus F : C\cap B_k[a]\ne
\emptyset\}$ is a local fragment. Since $C$ is special, $C\in X$
and therefore $C\in Q\setminus F$. \qed

Denote by $lc(H)$ the length of a longest induced cycle of the
graph $H$.

\begin{lemma}\label{mainchord} Let $G$ be a chordal graph with $kdim(G)\leq
3$. Let further there are no special cliques in $H$. Then
$lc(H)\leq 6$.
\end{lemma}

\pf Suppose contrary, i.e. let $a_1,\ldots,a_k$ form the induced
cycle $S \cong C_k$ in $H$, $k\geq 7$, $a_i a_{i+1}\in E(H)$,
indices are taken modulo $k$.

Since for every $a_i$ there are two nonadjacent neighbors in $H$,
then in every local fragment with center in $a_i$ it is covered by
at least 2 clusters. It implies $l_F(a_i)\leq 1$ for every
$i=1,\ldots,k$.

As $G$ is a chordal graph, there exist chords of the cycle $S$
covered by the fragment $F$. It is easy to see, that for every two
consecutive vertices $a_i,a_{i+1}$ of $S$ at least one of them
belongs to some chord of $S$ (indices are taken modulo $k$).
Indeed, let without lost of generality $a_i=a_k$, $a_{i+1} = a_1$.
If our statement is not true, then one can choose the chord
$a_pa_q$, $1 < p < q < k$ such, that $(p-1) + (k-q)$ is minimal.
But then $G(a_1,\ldots,a_p,a_q,\ldots,a_k)$ is a chordless cycle.

Assume without lost of generality, that one of chords of $S$
contains $a_1$. As $l_F(a_i)\leq 1$, for every vertex $a_i$ chords
incident to this vertex are covered by exactly one cluster of $F$.
It implies that there are no pairs of chords of the form
$\{a_ia_j, a_ia_{j+1}\}$, since in this case the vertices
$a_i,a_j,a_{j+1}$ are covered by one cluster of $F$ and thus the
edge $a_j a_{j+1}$ should be covered by $F$.

Let us show, that all chords of $S$ are covered by the cluster
$C_{chord} \supseteq \{a_1,a_3,\ldots,a_{k-1}\}$ (and thus $k$ is
even). Indeed, suppose that some chords of $S$ are covered by the
cluster $C\supseteq\{a_{i_1},\ldots,a_{i_r}\}$, $i_1 < i_2
<\ldots<i_r$, $i_1=1$, $C\ne C_{chord}$. Then there exist $1 \leq
p < q \leq r$ such, that $q-p \geq 3$. So, $G(a_{i_p},a_{i_p
+1},\ldots,a_{i_q-1},a_{i_q})$ is a cycle of length at least 4,
where without lost of generality $a_{i_p}$ belongs to some chord.
That chord should be covered by a cluster $C'\in F$, $C'\ne C$.
So, we have $l_F(a_{i_p}) \geq 2$, the contradiction.

In particular, this proposition implies that for any odd $i$ and
even $j$ such that $a_ia_j$ is not the edge of $S$, the vertices
$a_i$ and $a_j$ are nonadjacent in $G$ (otherwise $l_F(a_i)\geq
2$).

Let us denote by $C(v_1,\ldots,v_r)$ the clique containing
vertices $v_1,\ldots,v_r$.

Let $e=min\{ecc_H(a_1),5\}$. Since there is no special cliques in
$H$ there exists two different local fragments $F_e(a_1) \supseteq
\{C(a_1,a_2),C(a_1,a_k)\}$ and $F'_e(a_1)\supseteq
\{C'(a_1,a_2),C'(a_1,a_k)\}$, such that without lost of generality
$C(a_1,a_2)\setminus C'(a_1,a_2)\ne \emptyset$.

Let $v\in C(a_1,a_2)\setminus C'(a_1,a_2)$. Then $v\in
C'(a_1,a_k)$ and therefore $a_1v,a_2v,a_kv \in E(H)$.

\begin{figure}[h]\label{Figure1}
\begin{center}
\includegraphics*[width=40mm] {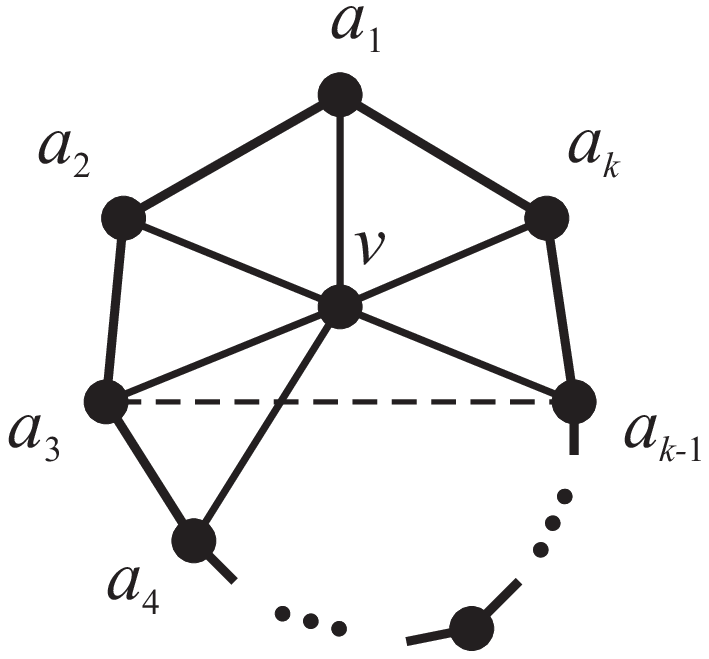}
\caption{}
\end{center}
\end{figure}

The vertices $a_2,v,a_k,a_{k-1},a_3$ form a cycle in $G$. It
should have at least 2 chords. Since $a_2a_{k-1},a_3a_k\not\in
E(G)$, there are edges $va_3,va_{k-1}\in E(G)$. The edges
$va_3,va_{k-1}$ are not covered by $F$ (otherwise $v\in C_{chord}$
and thus $\{a_1,v\}\in C_{chord}\cap C(a_1,a_2)$) and hence
$va_3,va_{k-1}\in E(H)$. It implies, that $v\in C(a_3,a_4)\in
F_e(a_1)$. So, $va_4\in E(H)$ (see Figure~1). Note, that since $k
\geq 7$ we have $a_4a_{k-1}\not\in E(H)$.

Let us remind, that in the local fragment $F'_e(a_1)$ the vertex
$v$ is covered by the cluster $C'(v,a_1,a_k)$. So, all other edges
of $H$ incident to $v$, should be covered by at most two clusters
of $F_e(a_1)$. But it is impossible, since the vertices $a_2$,
$a_4$, $a_{k-1}$ are pairwise nonadjacent. This contradiction
proves Lemma \ref{mainchord}. \qed

The considerations above suggest the following algorithm which
reduces the problem of recognition chordal graphs with krausz
dimension at most 3 to the same problem for graphs with bounded
maximum degree and maximum induced cycle length.

{\bf Algorithm 1}

{\bf Input:} chordal graph $G$.

{\bf Output:} One of the following:

1) graph $H$ with $\Delta(H)\leq 18$ and $lc(H)\leq 6$ such that
$kdim(G)\leq 3$ if and only if $kdim(H)\leq 3$;

2) the answer $"kdim(G) > 3"$.

{\bf begin}

\a $F:=\emptyset;$ $H:=G$; $isContinue := true$;

\a {\bf while} ($isContinue =$ {\bf true})

\b {\bf if} there exists a vertex $v\in V(H)$ such that $l_F(v) =
2$

\c $C:=N(v)\cup \{v\}$;

\c {\bf if} $C$ is a clique

\d  $F:=F\cup\{C\}$; {\bf continue} to the next iteration of the
cycle;

\c  {\bf else} the answer is "$kdim(G) > 3$"; {\bf stop;}

\b {\bf if} there exists a vertex $v\in V(H)$ with $deg(v)\geq 19$

\c  {\bf if} $v$ is contained in a clique $C$ with $|C|\geq 8$

\d  extend $C$ to a maximal clique; $F:=F\cup\{C\}$;

\d {\bf continue} to the next iteration of the cycle;

\c {\bf else} the answer is "$kdim(G) > 3$"; {\bf stop;}

\b For every non-isolated vertex $v\in V(H)$ generate all

\b $(v,e)$-local fragments, $e=min\{ecc_H(v),5\}$;

\b {\bf if} there exists a vertex $v\in V(H)$ such that there is

\c no $(v,e)$-local fragments

\c the answer is "$kdim(G) > 3$"; {\bf stop;}

\b {\bf if} there exists a special clique $C$

\c $F:=F\cup\{C\}$; {\bf continue} to the next iteration of the
cycle;

\b $isContinue := false$

\a {\bf endwhile};

\a add a pendant edge $vp_v$ to every vertex $v\in V(H)$ with
$l_F(v) = 1$;

{\bf end}.

\begin{tm}\cite{BT97}\label{boundtreewidth}
Let $lc(H)\leq s+2$, $\Delta(H)\leq \Delta$. Then $treewidth(H)
\leq \Delta(\Delta-1)^{s-1}$.

\end{tm}

\begin{tm} The problem $KDIM(3)$ is polynomially
solvable for chordal graph.
\end{tm}

\pf The correctness of algorithm 1 follows from the considerations
above. Let us show, that the Algorithm 1 is polynomial. Indeed,
the procedure of finding large clique which contains the fixed
vertex $v\in V(H)$ has the complexity $O(m)$. We start to generate
all possible $(v,e)$-local fragments for a vertex $v\in V(H)$ only
then $deg(v)\leq 18$. It implies $|B_e[v]|\leq const$ and thus the
complexity of this procedure is constant. The outer loop of the
algorithm 1 is performed at most $m$ times.

After performing the Algorithm 1 we obtained the graph $H$ with
bounded maximum degree and the length of a longest induced cycle.
By Theorem \ref{boundtreewidth} $H$ has bounded treewidth. For
such a graph the problem of determining its krausz dimension is
polynomially solvable \cite{HK97}.

\qed

\section{m-krausz dimension of graphs}\label{mkrausz}

We will start with proving the NP-hardness of the problem
$KDIM_m$. In order to make the proof more clear, we firstly will
prove, that $KDIM_m$ is NP-hard for general graphs, and then we
will use the developed construction to prove, that $KDIM_m$ is
NP-hard for $(1,2)$-colorable graphs.

\begin{tm}\label{kdimmNP}
The problem $KDIM_m$ is NP-hard for every fixed $m\geq 1$.
\end{tm}

\pf Let us reduce to the problem $KDIM_m$ the following special
case of the 3-dimensional matching problem (which we will call the
problem A):

Given: Non-intersecting sets $X$, $Y$, $Z$ , such that
$|X|=|Y|=|Z|=q$; $M\subseteq X\times Y\times Z$ , such that the
following condition holds:

(*) if $(a,b,w),(a,x,c),(y,b,c)\in M$, then $(a,b,c)\in M$.

The question:  Does $M$ contain a subset $M'\subseteq M$ ({\it
3-dimensional matching}) such, that $|M'|=q$ and every two
elements of $M'$ have not common coordinates?

It is known, that the problem A is NP-complete \cite{GJ}. Let $X$,
$Y$, $Z$, $M$, $|X|=|Y|=|Z|=q$, be the input of the problem A. Let
us reduce the problem A to  the problem of determining, if
$kdim_m(G)\leq 2q$. Construct the graph $G$ as follows:

\begin{equation}\label{vertreduce}
    V(G)=X\cup Y \cup Z \cup \{v,v_1,\ldots,v_q\};
\end{equation}

\begin{equation}\label{edgereduce}
    E(G)=\bigcup_{(a,b,c)\in M}\{ab,bc,ac\}\cup\{v_iv :
    i=1,\ldots,q\}\cup \{vd : d\in X \cup Y \cup Z\}
\end{equation}
(see Figure~2). Let us show that $M$ contains the 3-dimensional
matching $M'$ if and only if there exists a krausz
$(2q,m)$-partition of $G$.

\begin{figure}[h]\label{Figure2}
\begin{center}
\includegraphics*[width=50mm] {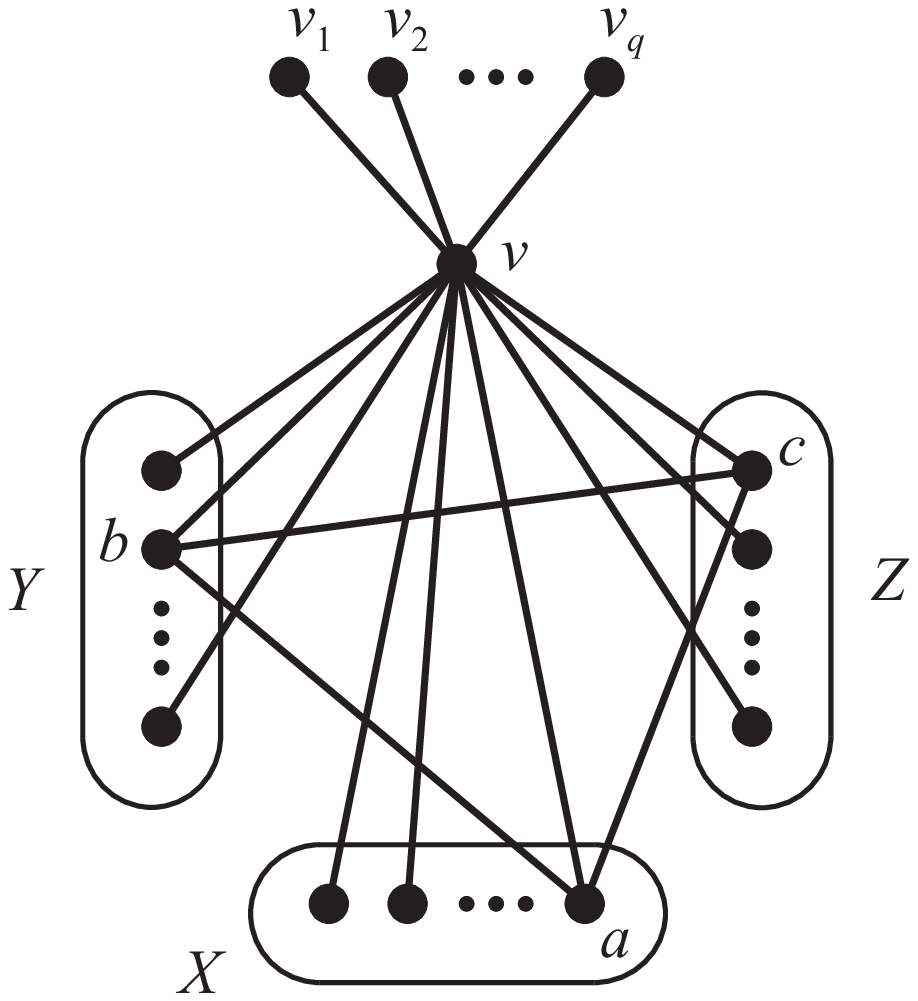}
\caption{}
\end{center}
\end{figure}

Suppose, that $M'=\{(a_i,b_i,c_i) : i=1,\ldots,q\}$ is the
3-dimensional matching. Let $Q_1 = \{\{v,a_i,b_i,c_i\} :
i=1,\ldots,q\}$, $Q_2 = \{\{v,v_i\} : i=1,\ldots,q\}$, $Q_3 =
\{\{z,t\} : zt\in E(G-(Q_1\cup Q_2))$. Then $Q=Q_1\cup Q_2\cup
Q_3$ is krausz $(2q,m)$-partition of $G$, since $deg(u)\leq 2q$
for every vertex $u\in V(G)\setminus\{v\}$ and the vertex $v$ is
covered by exactly $2q$ clusters of $Q$.

Let now $Q$ be krausz $(2q,m)$-partition of $G$. Denote by $Q(v)$
the set of clusters of $Q$, which contain the vertex $v$. Since
the vertices $v_i$, $i=1,\ldots,q$, have degree 1, there exist $q$
clusters from $Q(v)$ of the form $v v_i$, $i=1,\ldots,q$. Let
$C_1,\ldots,C_p$ be the remaining clusters from $Q(v)$, $p\leq q$.
Then $(C_1\cup\ldots\cup C_p)\setminus \{v\} = X\cup Y \cup Z$.
Since $X$, $Y$, $Z$ are stable sets of $G$, we have $|C_i|\leq 4$,
$i=1,\ldots,p$. As $|X\cup Y \cup Z| = 3q$, we have $p=q$,
$|C_i|=4$, $C_i\cap C_j = \{v\}$, $i,j=1,\ldots,p$, $i\ne j$.

Let $C_i = \{a_i,b_i,c_i,v : a_i\in X,\ b_i\in Y,\ c_i\in Z\}$,
$i=1,\ldots,q$. The property (*) implies, that $M' =
\{(a_i,b_i,c_i) : i=1,\ldots,q\}\subseteq M$ and, by the
consideration above, $M'$ is the 3-dimensional matching. \qed

\begin{cor}
The problem $KDIM_m$ is NP-hard in the class of $(1,2)$-colorable
graphs for every fixed $m\geq 1$.
\end{cor}

\pf Let us show, that the problem A could be reduced to the
problem $KDIM_m$ in the class of $(1,2)$-colorable graphs.

Let $G$ be the graph constructed in the proof of Theorem
\ref{kdimmNP}. Let us construct the graph $G'$ as follows:
$V(G')=V(G)\cup V_1'\cup V_2'$, where
\begin{equation}
    V_1' = \{w,w_1,\ldots,w_{2q}\};
\end{equation}
\begin{equation}
    V_2' = \{f_u : u\in V(G)\setminus (X\cup \{v_1,\ldots,v_q\})\};
\end{equation}
$E(G')=E(G)\cup E_1'\cup E_2'\cup E_3' \cup E_4'$, where
\begin{equation}
    E_1' = \{ww_i : i=1,\ldots,2q\};
\end{equation}
\begin{equation}
    E_2' = \{wx : x\in X\};
\end{equation}
\begin{equation}
    E_3' = \{x_1x_2 : x_1,x_2\in X, x_1\ne x_2\};
\end{equation}
\begin{equation}
    E_4' = \{uf_u : u\in V(G)\setminus (X\cup \{v_1,\ldots,v_q\}) \}
\end{equation}
(see Figure~3). The set $X\cup\{v\}$ is a clique, and the sets
$Y\cup \{f_z : z\in Z\}\cup \{v_1,\ldots,v_q, f_v,w\}$ and $Z\cup
\{f_y : y\in Y\}\cup \{w_1,\ldots,w_{2q}\}$ are stable sets of
$G'$. So, $G'$ is $(1,2)$-colorable graph.

\begin{figure}[h]\label{Figure3}
\begin{center}
\includegraphics*[width=75mm] {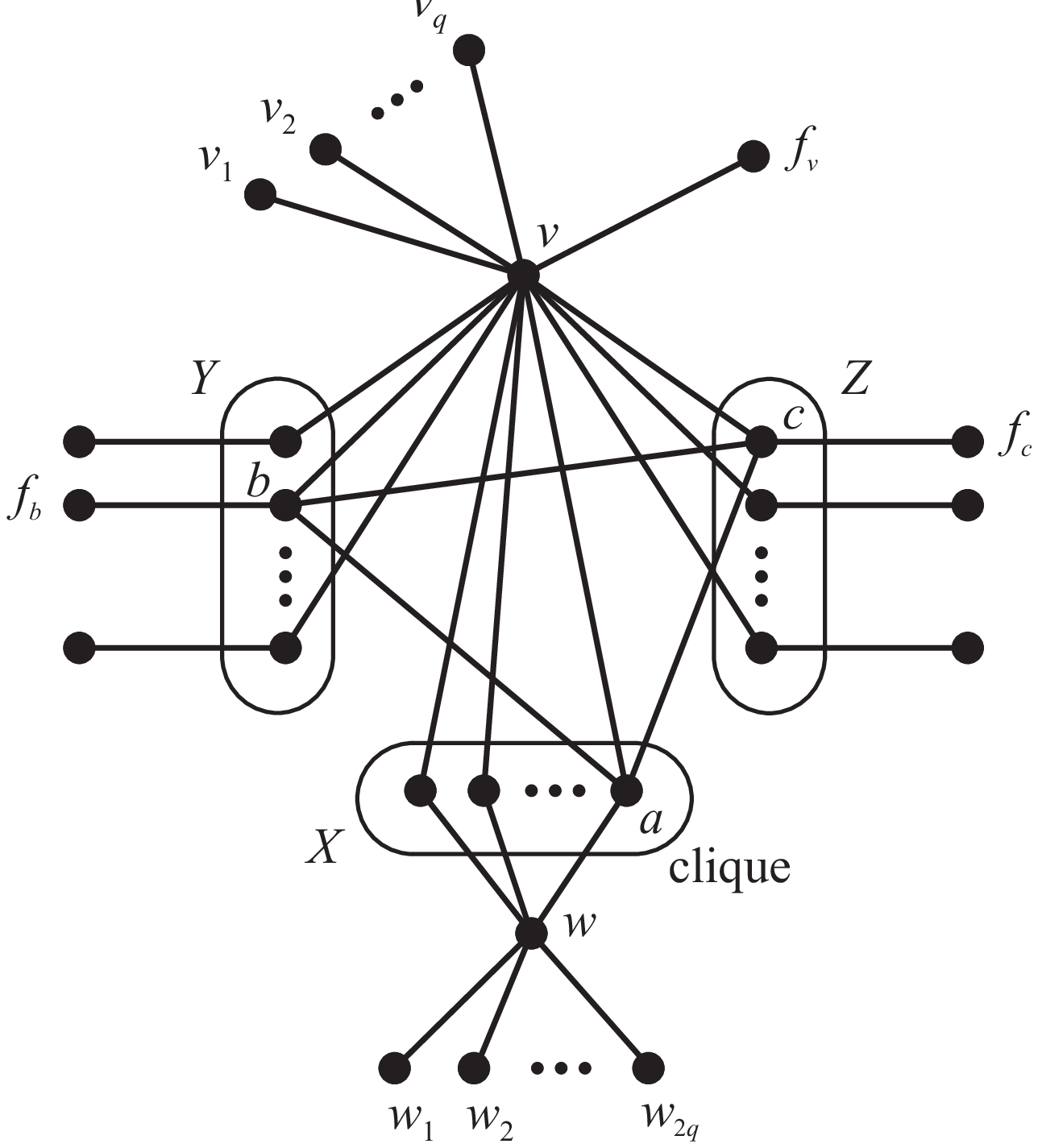}
\caption{}
\end{center}
\end{figure}

It is evident, that $Q$ is the krausz $(2q,m)$-partition of $G$ if
and only if
\begin{equation}
    Q\cup\{X\cup \{w\}\}\cup \{\{ww_i\} : i=1,\ldots,2q\}\cup \{\{u,f_u\}
    : u\in V(G)\setminus (X\cup \{v_1,\ldots,v_q\})\}
\end{equation}
is the krausz $(2q+1,m)$-partition of $G'$. \qed

Now we turn to the complexity of the recognition problem
$KDIM_m(k)$ in the class of $(\infty,1)$-polar graphs.

A maximal clique with at least $f(k,m)=m(k^2-k+1)+1$ vertices is
called a {\it $(k,m)$-large} clique.

In \cite{MSch08} the following two statements were proved. Since
they were published only in Russian in a journal, which is
difficult of access for a general reader, we repeat their proofs
here.

\begin{tm}\label{kmlarge} Any $(k,m)$-large
clique $C$ of a graph $G$ belongs to every krausz
$(k,m)$-partition of $G.$
\end{tm}

\pf Let $A$ be a krausz $(k,m)$-partition of graph $G$,
$A_1,A_2,\ldots,A_t$ be those clusters of $A$ which have common
vertices with $C$. Assume that $C\not\in A$. Then the family $B =
(B_1,B_2,\ldots,B_t)$, where $B_i = A_i\cap C$, is a krausz
$(k,m)$-partition of the graph $G(C)$, and (by maximality of $C$)
$B_i\ne C$ for every $i=1,2,\ldots,t$.

Let us show, that $|B_i|\le mk$ for any $i=1,2,\ldots,t$. Consider
a cluster of $B$, say $B_1$, and a vertex $u \in C\setminus B_1$.
No edge of the form $ux$, where $x \in B_1 $, is contained in
$B_1$. Moreover, each cluster of $B$ different from $B_1$ contains
at most $m$ of such edges (by the definition of krausz
$(k,m)$-partition). Taking into account that the vertex $u$
belongs to at most $k$ clusters of $B$, we obtain the inequality
$|B_1|\le mk$.

Now we will prove that if $B_i\setminus B_j\ne\emptyset$ for some
clusters $B_j\in B$, then $|B_j\setminus B_i|\le m(k - 1)$.
Consider a vertex $u\in B_i\setminus B_j$. Any edge of the form
$ux$, where $x\in B_j\setminus B_i$ (if such one exists) is
contained neither in $B_i$, nor in $B_j$. Besides, no cluster of
$B$ contains more than $m$ of such edges by definition of krausz
$(k,m)$-partition. Taking into account that $u$ belongs to at most
$k-1$ clusters of $B$ different from $B_i$, we obtain the
inequality $|B_j \setminus B_i|\le m(k-1)$.

Consider an arbitrary vertex $v$ of the clique $C$. Let, without
loss of generality, it belongs to the clusters
$B_1,B_2,\ldots,B_s$ of $B$, $s\le t$. We show that $|B_1\cup
B_2\cup\ldots\cup B_s|\le mk+(s-1)m(k-1)$. The following equality
is obvious
\begin{equation}
\label{eq1} |B_1\cup B_2\cup\ldots\cup B_s| = |B_1| +
|B_2\setminus B_1| + |B_3\setminus (B_1\cup B_2)| + \ldots +
|B_s\setminus (B_1\cup B_2\ldots\cup B_{s - 1} )|.
\end{equation}
If $B_1\setminus B_2\ne\emptyset$, $(B_1\cup B_2)\setminus
B_3\ne\emptyset$, \ldots, $(B_1\cup B_2\cup\ldots\cup
B_{s-1})\setminus B_s\ne\emptyset$, then by proved above each term
in the right part of the equality (\ref{eq1}), starting from the
second, does not exceed $m(k-1)$. Hence we have $|B_1\cup
B_2\cup\ldots\cup B_s|\le mk+(s-1)m(k-1)$. Let, on the contrary,
$i\in \{2,\ldots,s\}$ is the maximal number such, that
$(B_1\cup\ldots\cup B_{i-1})\setminus B_i =\emptyset$. Then
$B_1\subseteq B_i$, $B_2\subseteq B_i$, \ldots, $B_{i -
1}\subseteq B_i$, and the sum of the first $i$ terms in the right
part of (\ref{eq1}) is equal to $|B_1\cup B_2\cup\ldots\cup
B_i|=|B_i|\le mk$. Each of the other terms does not exceed
$m(k-1)$ by the maximality of $i$. Hence
$$|B_1\cup B_2\cup\ldots\cup B_s|\le mk+(s-i)m(k - 1) < mk+(s-1)m(k - 1).$$

So, in any case we obtain that the inequality $|B_1\cup
B_2\cup\ldots\cup B_s|\le mk+(s-1)m(k - 1)$ holds. Taking into
account that $C = B_1\cup B_2\cup\ldots\cup B_s$ and $s\le k$, we
have
$$|C|\le mk+(k-1)m(k-1) = m(k^2-k+1) < f(k,m).$$
The obtained contradiction proves the lemma.   \qed

\begin{tm}\label{kdimmsplit}
There exists a finite set $\mathcal{F}_0$ of forbidden induced
subgraphs such that a split graph $G$ belongs to the class $L_k^m$
if and only if no induced subgraph of $G$ is isomorphic to an
element of $\mathcal{F}_0$.
\end{tm}
\pf Denote by $R_p$ the graph obtained from the complete graph
$H\cong K_{f(k,m)}$ by adding a new vertex and connecting it with
exactly $p$ vertices of $H$. Put $\mathcal{F}_0 = \{R_p : km+1\le
p\le f(k,m)-1 \}\cup\{K_{1,k+1}\}$. Using Theorem \ref{kmlarge}
one can immediately verify that no graph from $\mathcal{F}_0$
belongs to $L_k^m$.

Let, without loss of generality, $G$ be connected graph, and $V(G)
= C\cup S$ be a bipartition of $V(G)$ into clique $C$ and stable
set $S$ such, that $C$ is a maximal clique. Let also no induced
subgraph of $G$ be isomorphic to an element of $\mathcal{F}_0$.
Put $S=\{v_1,\ldots,v_s\}$. Consider two cases:

1) $|C| > (km-1)k+1$.

In this case we have
$$|C| \ge (km-1)k+2 = mk^2-(k-1)+1 \ge mk^2-m(k-1)+1 = f(k,m).$$
Then, since no induced subgraph of $G$ is isomorphic to a graph
$R_p$, $km+1\le p\le f(k,m)-1$, we have $\deg(v_i)\le km$ for any
$i=1,2,\ldots,s$. Since $G$ contains no induced $K_{1,k+1}$, we
have $|N(u)\cap S|\le k$ for any vertex $u$ from $C$. Moreover, we
prove that for any vertex $u$ from $C$ the inequality $|N(u)\cap
S|\le k-1$ holds. Assume this is not true. Let, without loss of
generality, some vertex $u$ from $C$ be adjacent to the vertices
$v_1,\ldots,v_k$ from $S$, $k\le s$. Since $deg(v_i)\le km$,
$i=1,2,\ldots,k$, and $u\in\bigcap\limits_{i=1}^k N(v_i)$, then
$|\bigcup\limits_{i=1}^k {N(v_i)}| \le \sum\limits_{i=1}^k
(deg(v_i)-1)+1 \le (km-1)k + 1 < \varphi(G)$. Hence, there exists
a vertex $u'$ from $C$, which is not adjacent to any vertex from
$v_1,\ldots,v_k$. But then $G(u,u',v_1,\ldots,v_k)\cong
K_{1,k+1}$, a contradiction.

Now we can construct a krausz $(k,m)$-partition of $G$. Since
$deg(v_i)\le km$ for any $i=1,2,\ldots,s$, then there exists a
partition $N(v_i) = C_{i_1}\cup\ldots\cup C_{i_{s_i}}$, where
$C_{i_j}\cap C_{i_l} =\emptyset$, $j,l\in\{1,\ldots,s_i\}$, $j\ne
l$, $|C_{i_j}|\le m$, $s_i\le k$. Obviously, the list of cliques
$\{C_{i_j}\cup \{v_i\} : i = \overline{1,s}, j =
\overline{1,s_i}\}$ together with the clique $C$ is a krausz
$(k,m)$-partition of graph $G$.

2) $|C| \le (km-1)k+1$.

Since $G$ contains no induced $K_{1,k+1}$, we have $|N(u)\cap
S|\le k$ for any vertex $u$ from $C$. Therefore, as $G$ is
connected,
$$
|G| = |C|+|S| \le |C|+\sum\limits_{u \in C} |N(u)\cap S| \le
((km-1)k+1)+((km-1)k+1)k = ((km-1)k+1)(k+1),
$$
i.~e. the order of graph $G$ is bounded above by a value,
depending on $k$ and $m$. Add to the list $\mathcal {F}_0$ all
such split graphs $H$, that $H\not\in L_k^m$ and $|H| \le
((km-1)k+1)(k+1)$.

Obviously, the constructed in the cases 1) and 2) finite list
$\mathcal {F}_0$ is a required list of forbidden induced
subgraphs. \qed

Since $K_{1,k+1}\not\in L_k^m$, the heredity of $L_k^m$
immediately implies
\begin{lemma}\label{kdimmbipartite}
A bipartite graph $G$ belongs to the class $L_k^m$ if and only if
no induced subgraph of $G$ is isomorphic to $K_{1,k+1}$.
\end{lemma}

\begin{tm}\label{kdimmpolar}
There exists a finite set $\mathcal{F}$ of forbidden induced
subgraphs such that an $(\infty,1)$-polar graph $G$ belongs to the
class $L_k^m$ if and only if no induced subgraph of $G$ is
isomorphic to an element of $\mathcal{F}$.
\end{tm}

\pf Without loss of generality we can suppose that
$(\infty,1)$-polar graph $G$ is connected. Let $G$ have
bipartition $(A,B)$; $A_i$, $i=1,2,\ldots,t$, be the vertex sets
of connected components of $\overline{G}(A)$; $\mathcal{F}_0$ be
the set of split graphs from Theorem \ref{kdimmsplit}. Denote by
$\mathcal{F}_1$ the set of $(\infty,1)$-polar graphs which do not
belong to the class $L_k^m$ and have order at most
$(k+1)k(f(k,m)-1)$.

Put $\mathcal{F}=\mathcal{F}_0\cup \mathcal{F}_1\cup \{K_{1,k+1},
K_{f(k,m)+1}-e\}$, where $K_{f(k,m)+1}-e$ is the graph obtained
from the complete graph $K_{f(k,m)+1}$ after deleting an edge. The
set $\mathcal{F}$ is finite, since $\mathcal{F}_0$ and
$\mathcal{F}_1$ are finite. According to Theorem \ref{kmlarge},
there is no krausz $(k,m)$-partition for $K_{f(k,m)+1}-e$.
Therefore $K_{f(k,m)+1}-e\not\in L_k^m$. Thus, $\mathcal{F}\cap
L_k^m=\emptyset$. The necessity of the statement follows from the
heredity of the class $L_k^m$.

Now let $G$ contain no induced subgraph isomorphic to an element
from $\mathcal{F}$. If $G(A)$ is complete, then $G$ is split graph
and by Theorem \ref{kdimmsplit} $G\in L_k^m$. If $G(A)$ is empty,
then $G$ is bipartite graph and by Lemma \ref{kdimmbipartite}
$G\in L_k^m$.

Now suppose that $G(A)$ is neither complete nor bipartite graph.
Then $2\le t\le |A|-1$. Since $K_{1,k+1}\in \mathcal{F}$, then
$|A_i|\le k$ for any $i=1,2,\ldots,t$. Now we will prove that
since $K_{f(k,m)+1}-e\in \mathcal{F}$, then $t\le f(k,m)-1$. Let,
to the contrary, $t\ge f(k,m)$. As $G(A)$ is not complete graph,
there exists an index $i_0\in \{1,2,\ldots,t\}$ such that
$|A_{i_0}|\ge 2$. Consider the set
$S=\{a_1,a_2,\ldots,a_{i_0-1},a'_{i_0},a''_{i_0},a_{i_0+1},\ldots,a_t\}$,
where $a_i\in A_i$ for any $i\in \{1,2,\ldots,t\}\setminus
\{i_0\}$ and $a'_{i_0},a''_{i_0}\in A_{i_0}$. Then $G(S)$ contains
$K_{f(k,m)+1}-e$ as induced subgraph, a contradiction. Therefore
$$|A|\le \sum_{i=1}^t |A_i|\le k(f(k,m)-1).$$
Since $|N(a)\cap B|\le k$ for any vertex $a\in A$ and $G$ is
connected, we have
$$|G|\le |A|+|B|\le |A|+ \sum_{a\in A} |N(a)\cap B|\le
k(f(k,m)-1)+k^2(f(k,m)-1)=(k+1)k(f(k,m)-1).$$ It follows from the
inclusion $\mathcal{F}_1\subseteq \mathcal{F}$ that $G\in L_k^m$.
\qed

\begin{cor}\label{polinominpolar}
The problem $KDIM_m(k)$ is polynomially solvable in the class of
$(\infty,1)$-polar graphs for every fixed $k,m\geq 1$.
\end{cor}


\end{document}